\documentclass[12pt,reqno]{amsart}
\usepackage{amssymb}
\usepackage{longtable}
\theoremstyle{plain}
\newtheorem{theorem}{Theorem}
\newtheorem{lemma}{Lemma}
\newtheorem{corollary}{Corollary}[section]
\newtheorem{proposition}{Proposition}[section]
\theoremstyle{proof}
\theoremstyle{definition}

\newtheorem{definition}{Definition}[section]

\theoremstyle{remark}
\newtheorem{remark}{Remark}
\theoremstyle{lamma}

\numberwithin{equation}{section}
\numberwithin{lemma}{section}
\numberwithin{theorem}{section}

\usepackage{amssymb, amsmath, amsthm}
\usepackage[breaklinks]{hyperref}
\theoremstyle{thmrm}

\usepackage{graphicx}
\begin{document}
 \title[]{Irreducible integrable Modules for the full toroidal Lie algebras co-ordinated by Rational Quantum Torus}
\address{Harish-Chandra Research Institute,  A CI of Homi Bhabha National
Institute, Chhatnag Road, Jhunsi, Prayagraj - 211019}
\email{santanutantubay@hri.res.in, batra@hri.res.in}
 
\author{Santanu Tantubay, Punita Batra}
\maketitle
\begin{abstract}
Let $\mathbb{C}_q$ be the non-commutative Laurent polynomial ring associated with a $(n+1)\times (n+1)$ rational quantum matrix . Let $\mathfrak{sl}_d(\mathbb{C}_q)\oplus HC_1(\mathbb{C}_q)$ be the universal central extension of Lie subalgebra $\mathfrak{sl}_d(\mathbb{C}_q)$ of $\mathfrak{gl}_d(\mathbb{C}_q)$ . Now let us take the Lie algebra $\tau=\mathfrak{gl}_d(\mathbb{C}_q)\oplus HC_1(\mathbb{C}_q)$. Let $Der(\mathbb{C}_q)$ be the Lie algebra of all derivations of $\mathbb{C}_q$. Now we consider the Lie algebra $\tilde{\tau}=\tau\rtimes Der(\mathbb{C}_q)$, called as full toroidal Lie algebra co-ordinated by rational quantum tori. In this paper we get a classification of irreducible integrable modules with finite dimensional weight spaces for $\tilde{\tau}$ with nonzero central action on the modules.
\end{abstract}
\section{introduction}
EALAs of type $A_{d-1}$ are tied up with the Lie algebra $\mathfrak{gl}_d(\mathbb{C}_q)$ corresponding to some rational quantum matrix $q$. Co-ordinate ring of such algebra is non-commutative analogue of Laurent polynomial algebras. Non zero level integrable modules with finite dimensional weight spaces of such algebras were classified in \cite{[2]}. It is well known that the derivation algebra of non-commutative Laurent polynomial ring naturally acts on $\mathfrak{gl}_d(\mathbb{C}_q)$. Note that when all $q_{ij}=1$, $\mathbb{C}_q$ becomes the commutative Laurent polynomial ring $A$ in $(n+1)$ variables $t_0,\dots,t_n$. In \cite{[4]}, S. Eswara Rao classified irreducible integrable modules for toroidal Lie algebra with commutative co-ordinate algebra, where each weight space is of finite dimension. Then in \cite{[6]}, S. Eswara Rao and C. Jiang consider the full toroidal Lie algebra with commutative co-ordinate algebra and classified such module for this Lie algebra. In the case of full toroidal Lie algebra, they identified the highest weight space with an irreducible $\mathbb{Z}^{(n+1)}$-graded  $A\rtimes Der(A)$-module. In \cite{[19]}, S. Eswara Rao classified irreducible modules for $A\rtimes Der (A)$ with associative $A$-action on the modules. In \cite{[3]}, S. Eswara Rao, Punita Batra, Sachin. S. Sharma classified irreducible $\mathbb{Z}^{(n+1)}$-graded  modules for $\mathbb{C}_q\rtimes Der(\mathbb{C}_q)$ with associative and anti-associative action of $\mathbb{C}_q^1$ and $\mathbb{C}_q^2$ respectively (see \cite{[3]} for the defination of $\mathbb{C}_q^i$ for $i=1,2$) . In this paper, we consider the full toroidal Lie algebra co-ordinated by rational quantum torus $\mathbb{C}_q$ and classify irreducible integrable modules for this Lie algebra with finite dimensional weight spaces. In \cite{[20]}, C. Jiang and D. Meng proved that classification of irreducible integrable modules for full toroidal Lie algebra with finite dimensional weight spaces reduces to the problem of classification of irreducible modules for $A\rtimes Der(A)$ but in our case we do not have this type of correspondence, so here we can not identify the highest weight space with an irreducible $\mathbb{Z}^{(n+1)}$-graded $\mathbb{C}_q\rtimes Der(\mathbb{C}_q)$ module. In this case non-commutativity of the co-ordinate algebra is troublesome. So in order to classify such modules we need some different approach. In this case we are taking the help of \cite{[7]} and able to settle the problem. Then our approach is similar to \cite{[11]}. \\
Now we summarize the content of this article. In Section \ref{section 2}, we define a rational quantum matrix $q$ and the corresponding quantum torus $\mathbb{C}_q$ associated to the matrix $q$. After recalling some results about $\mathbb{C}_q$, we define the Lie algebra $\mathfrak{gl}_d(\mathbb{C}_q)$ and its universal central extension $\tau$ and add derivation algebra of $\mathbb{C}_q$ with $\tau$. Now the resultant Lie algebra $\tilde{\tau}$ has a root space decomposition. In Section \ref{section 3}, we define integrable modules for $\tilde{\tau}$. Then we take a triangular decomposition of $\tilde{\tau}$ and showed that any non zero level irreducible integrable module with finite dimensional weight spaces is either a highest weight module or a lowest weight module. We prove the highest weight space of a highest weight module is irreducible for the zero-th component of triangular decomposition of $\tilde{\tau}$. In Section \ref{section 4}, we showed that a large part of central extension acts trivially on highest weight space. Rest of the part which sits inside zero-th component of triangular decomposition acts associatively on the highest weight space. In Section \ref{section 5}, we consider a Lie subalgebra $\mathfrak{L}$, of zero-th component of triangular decomposition and take a submodule $W$ of the highest weight space for $\mathfrak{L}$ and prove that the quotient module is completly reducible finite dimensional for $\mathfrak{L}$. In Section \ref{section 6}, we break this Lie algebra into two parts. To identify the first part of the Lie algebra we first give an isomorphism of $\mathbb{C}_{q^\prime}$ with a twisted multi-loop algebra. Using this isomorphism we can identify finite dimensional simple modules for first part of the Lie algebra. Then we prove that every irreducible components of the quotient module is tensor product of irreducible modules for $\oplus \mathfrak{gl}_N$ and $\mathfrak{gl}_n$. Then we prove our main result Theorem \ref{theorem 6.1}.  
\section{notations and preliminaries}\label{section 2}
In this paper all the vector spaces, algebras, tensor products are over the field of complex numbers $\mathbb{C}$. Let $\mathbb{Z},\mathbb{N}$ denote the set of integers and natural numbers respectively. For any Lie algebra $L$, let $U(L)$ denote the universal enveloping algebra of $L$ and $\mathbb{C}^n$ be the $n$ copies of complex field $\mathbb{C}$. Let $\{e_i:1\leq i\leq n\}$ be the standard basis of $\mathbb{C}^n$ and $(,)$ be the standard form on $\mathbb{C}^n$.
Let us fix  two positive integers $d,n\geq 2$. Let $q=(q_{i,j})$ be a $(n+1)\times (n+1)$ matrix such that $q_{i,i}=1$, $q_{i,j}=q_{j,i}^{-1}$ and $q_{i,j}$ is a root of unity, where $0\leq i,j\leq n$. Let $\mathbb{C}_q$ be the Laurent polynomial ring in $(n+1)$ non-commutating variables $t_0,t_1,\dots ,t_n$ with the conditions $t_it_j =q_{i,j}t_jt_i$ for $0\leq i,j\leq n$. Clearly we can see that $\mathbb{C}_q$ is $\mathbb{Z}^{n+1}$-graded with each graded component is one dimensional. For any $a=(a_0,\dots, a_n) \in \mathbb{Z}^{n+1}$, let $t^a=t_0^{a_0}t_1^{a_1}\cdots t_n^{a_n} \in \mathbb{C}_q$.\\ Now let us define the following maps $\sigma, f :\mathbb{Z}^{n+1}\times \mathbb{Z}^{n+1} \rightarrow \mathbb{C}^*$ by
\begin{equation}
\sigma (a,b)= \prod _{0\leq i\leq j\leq n} q_{j,i}^{a_jb_i}
\end{equation}
\begin{equation}
f(a,b)=\sigma (a,b) \sigma (b,a)^{-1}
\end{equation} 
Then one has the following results for any $a,b,c \in \mathbb{Z}^{n+1}, k \in \mathbb{Z}$ :
\begin{enumerate}
\item $f(a,b)=f(b,a)^{-1}$,
\item $f(ka,a)=f(a,ka)=1$,
\item $f(a+b,c)=f(a,c)f(b,c)$,
\item $f(a,b+c)=f(a,b)f(a,c)$,
\item $\sigma (a,b+c)=\sigma (a,b)\sigma (a,c)$,
\item $t^at^b=\sigma(a,b)t^{a+b}$, $t^at^b=f(a,b)t^bt^a$.
\end{enumerate}
The radical of $f$ is defined by 
\begin{center}
$\text{rad} \;f=\{a \in \mathbb{Z}^{n+1}|\;f(a,b)=1, \forall\; b \in \mathbb{Z}^{n+1}\}$.
\end{center}
 It is easy to see that rad $f$ is a subgroup of $\mathbb{Z}^{n+1}$ and $m \in \text{rad}\; f$ iff $f(r,s)=1,\; \forall\; r,s \in \mathbb{Z}^{n+1}$ with $r+s=m$. Now $\mathbb{C}_q$ being an associative algebra, it has a natural Lie algebra structure. Now we recall a proposition from \cite{[1]}.
\begin{proposition}{(Proposition 2.44,[ \cite{[1]}])}
\begin{enumerate}
\item The center $Z(\mathbb{C}_q)$ of $\mathbb{C}_q$ has a basis consisting of monomials $t^a$, $a \in \text{rad}\; f $.
\item The Lie subalgebra $[\mathbb{C}_q,\mathbb{C}_q]$ of $\mathbb{C}_q$ has a basis consisting of monomial $t^a, a\in \mathbb{Z}^{n+1} \setminus \text{rad}\; f$.
\item $\mathbb{C}_q=[\mathbb{C}_q,\mathbb{C}_q]\oplus Z(\mathbb{C}_q)$.
\end{enumerate}
\end{proposition} 
Let $M_d(\mathbb{C})$ be the associative algebra of $d\times d$ matrices. Now consider the associative algebra $M_d(\mathbb{C}_q)=M_d(\mathbb{C})\otimes \mathbb{C}_q$, denote the corresponding Lie algebra by $\mathfrak{gl}_d(\mathbb{C}_q)$. For any $x \in M_d(\mathbb{C})$ and $t^a\in \mathbb{C}_q$, let $x(a)=x \otimes t^a$. Then the Lie bracket is given by 
\begin{center}
$[x(a),y(b)]_0=(\sigma(a,b)xy-\sigma (b,a)yx)(a+b)$.
\end{center}
Now we recall the universal central extension of $\mathfrak{gl}_d(\mathbb{C}_q)$ from \cite{[1]}. Let $J$ be the linear span of 
\begin{center}
$x\otimes y+y\otimes x$,   $xy\otimes z+yz\otimes x+zx\otimes y$
\end{center}
inside $\mathbb{C}_q\otimes \mathbb{C}_q$ for all $x,y,z \in \mathbb{C}_q$. Let $\langle x,y\rangle_0$ denote the element $x\otimes y+J$ in $\mathbb{C}_q \otimes \mathbb{C}_q/J$. Define
\begin{center}
\[HC_1(\mathbb{C}_q)= \Big\{\sum_i \langle x_i,y_i\rangle _0|\sum_i[x_i,y_i]=0\Big\}.\]
\end{center}
Let $\tau=\mathfrak{gl}_d(\mathbb{C}_q)\oplus HC_1(\mathbb{C}_q)$, where the Lie brackets are given by 
\begin{center}
$[x(a),y(b)]=[x(a),y(b)]_0+\text{Tr}(xy)\;\langle t^a,t^b\rangle_0\;\delta_{a+b,\text{rad}\; f},$
\end{center}
 where 
 \begin{center}
 \[ \delta_ {a,  \text{rad}\;f}=\begin{cases} 1 & \text{if}\; a \in \text{rad}\; f, \\
 0 & \text{if}\; a \notin \text{rad} \;f. \end{cases}\]
 \end{center}
From now on we write $\langle t^a,t^b\rangle$ for $\langle t^a,t^b\rangle_0\;\delta_{a+b, \text{rad}\; f}$.\\
We recall the following Lemma from \cite{[1]} about the properties of elements in $HC_1(\mathbb{C}_q)$.
\begin{lemma}\label{lemma 2.1}
\begin{enumerate}
\item $\langle 1,t^a\rangle=0$ for all $a \in \mathbb{Z}^{n+1}$.
\item $\langle t^b,t^a(t^b)^{-1}\rangle=\sum_{i=0}^{n}b_i\;\langle t_i,t^at_i^{-1}\rangle$.
\item $\langle t^a,t^b\rangle=\sigma (a,b)\sum_{i=0}^{n}a_i\;\langle t_i,t^{a+b}t_i^{-1}\rangle.$
\item $(t^b)^{-1}=\sigma(b,b)t^{-b}.$
\item $dim\; (HC_1(\mathbb{C}_q)_a=\begin{cases} 0& \text{if}\; a \notin \text{rad}\;f\\
n &\text{if}\; a \in \text{rad}\;f \setminus\{0\}\\n+1& \text{if}\; a=0\end{cases}$
\end{enumerate}
\end{lemma}
Let $Der(\mathbb{C}_q)$ be the space of all derivations of $\mathbb{C}_q$. Then we have the lemma
\begin{lemma}(See \cite{[1]}, Lemma 2.48)
\begin{enumerate}
\item Der$(\mathbb{C}_q)=\oplus_{a\in \mathbb{Z}^{n+1}}(\text{Der}(\mathbb{C}_q))_a$.
\item \[ \text{Der}(\mathbb{C}_q)_a=\begin{cases} \mathbb{C} \text{ad}\; t^a & \text{if}\;a \notin \text{rad}\;(f)\\ \bigoplus_{i=0}^{n}\mathbb{C}t^a  \partial_i & \text{if}\; a \in \text{rad}\;(f).\end{cases}\]
\end{enumerate}
\end{lemma} 
We denote $D(u,r)=\sum_{i=0}^{n}u_it^r\partial_i$ for $u\in \mathbb{C}^{n+1},\;r\in \text{rad}\;f$. The space Der($\mathbb{C}_q$) is a Lie algebra with the following brackets:
\begin{enumerate}
\item $[\text{ad}\;t^r.\text{ad}\;t^s]=(\sigma (r,s)-\sigma (s,r))\;\text{ad}\; t^{s+r},\: \forall\; r,s \notin \text{rad} \;f$;
\item $[D(u,r),\text{ad} \;t^s]=(u,s)\sigma (r,s)\;\text{ad}\; t^{r+s},\: \forall\; r \in \text{rad}\;(f), s \notin \text{rad} (f), u \in \mathbb{C}^{n+1}$;
\item $[D(u,r),D(u^{\prime},r^{\prime})]=D(w, r+r^{\prime}),\; \forall\; r,r^{\prime}\in \text{rad}\;(f), u,u^{\prime}\in \mathbb{C}^{n+1}$ and where $w=\sigma(r,r^\prime)((u,r^{\prime})-(u^{\prime},r)).$
\end{enumerate}
Now Der$(\mathbb{C}_q)$ naturally acts on $\mathfrak{gl}_d(\mathbb{C}_q)$. 
We recall from \cite{[1]}, that there is a natural action of $Der(\mathbb{C})_q$ on $HC_1(\mathbb{C}_q)$. The action is given as follows:
\[D(u,r). \langle t^a,t^b\rangle=(u,a)\sigma(r,a)\langle t^{a+r},t^b\rangle+(u,b)\sigma(r,b)\langle t^a,t^{b+r}\rangle;\]
\[ad\;t^s. \langle t^a,t^b\rangle=0\;\text{(Equation 2.57 of \cite{[1]})}.\]  Now we define the Lie algebra $\tilde{\tau}=\tau \rtimes \text{Der}(\mathbb{C}_q)$ with the brackets
\begin{enumerate}
\item $[D(u,r),x(a)]=(u,a)\sigma (r,a)\;x(r+a), \; \forall\; r \in \text{rad}\;(f),a \in \mathbb{Z}^{n+1},u\in \mathbb{C}^{n+1}$;
\item $[ad\; t^s,x(a)]=(\sigma (s,a)-\sigma (a,s))x(a+s),\; \forall\; s \notin \text{rad}(f),a \in \mathbb{Z}^{n+1};$
\item $[ad\;t^s,ad\;t^r]=(\sigma(s,r)-\sigma(r,s))ad\; t^{s+r},\; \forall\; r,s \notin \text{rad}(f);$
\item $[D(u,r),ad \; t^s]=(u,s)\sigma(r,s)ad\; t^{r+s},\; \forall\; r \in \text{rad}(f),s \notin \text{rad}(f),u \in \mathbb{C}^{n+1}$;
\item $[D(u,r),D(u^{\prime},r^{\prime})]=D(w,r+r^{\prime}),\;\forall\; r,r^{\prime} \in \text{rad}(f),u,u^{\prime}\in \mathbb{C}^{n+1}$ and where $w=\sigma(r,r^{\prime})((u,r^{\prime})u^{\prime}-(u^{\prime},r)u)$;
\item $[D(u,r),\langle t^a,t^b\rangle]=(u,a)\sigma(r,a)\langle t^{a+r},t^b\rangle+(u,b)\sigma(r,b)\langle t^a,t^{b+r}\rangle$, where $a,b \in \mathbb{Z}^{n+1}$ and $u\in \mathbb{C}^{n+1},\;r\in \text{rad}\;f$;
\item $[ad\; t^s,\langle t^a,t^b\rangle]=0$, where $a,b \in \mathbb{Z}^{n+1}$ and $s\notin \text{rad}\;f$.
\end{enumerate}
 Let us denote the zero degree elements of $HC_1(\mathbb{C}_q)$ by  $C_i=\langle t_i,t_i^{-1}\rangle$ for $i=0,\dots, n$. Let $\mathfrak{h}$ be the Cartan subalgebra of $\mathfrak{sl}_d(\mathbb{C})$ spanned by the elements $E_{ii}-E_{i+1,i+1}$ for $1\leq i\leq d-1$ and assume $\mathfrak{sl}_d(\mathbb{C})=(\mathfrak{sl}_d(\mathbb{C}))_-\oplus \mathfrak{h}\oplus (\mathfrak{sl}_d(\mathbb{C}))_+$ be the triangular decomposition for the simple Lie algebra $\mathfrak{sl}_d(\mathbb{C})$ with the cartan subalgebra $\mathfrak{h}$. Let $I$ be the identity matrix of $\mathfrak{gl}_d(\mathbb{C})$ and let $D$ be the $\mathbb{C}$ linear span of $\{\partial_i|0\leq i\leq n\}$. Let $\mathfrak{H}=\mathfrak{h}\oplus \mathbb{C}I\oplus \sum_{i=0}^n\mathbb{C}C_i \oplus D$, which will be a cartan subalgebra of $\tilde{\tau}$. Define $\delta_i,\;w_i,\;I^*,\;0\leq i\leq n$, such that $\delta_i(\mathfrak{h})=\delta_i(C_j)=\delta_i(I)=0,\; \delta_i(\partial_j)=\delta_{ij}$, $w_i(\mathfrak{h})=w_i(\partial_j)=w_i(I)=0,\;w_i(C_j)=\delta_{ij}$ and $I^*(\mathfrak{h})=I^*(C_i)=I^*(\partial_i)=0,\;I^*(I)=1$. Let $\alpha_1,\dots, \alpha_{d-1}$ be the standard simple roots of $\mathfrak{h}$. Let $\bigtriangleup _0$ be the root system for $\mathfrak{sl}_d(\mathbb{C})$. Let $\gamma =\alpha +\delta _m$ for $\alpha \in \bigtriangleup _0$. Then $\gamma$ is called real root. Let $\bigtriangleup ^{re}$ be the set of real roots and let $\bigtriangleup=\{\alpha+\delta_m,\delta_m|\alpha\in \bigtriangleup ^{re},\delta _m \;\text{is a null root}\}$. Let $\mathfrak{sl}_d(\mathbb{C})=\oplus_{\alpha \in \bigtriangleup_0 \cup \{0\}}(\mathfrak{sl}_d(\mathbb{C}))_{\alpha}$ be the root space decomposition.\\
 Define $\tilde{\tau}_{\alpha+\delta_m}=(\mathfrak{sl}_d(\mathbb{C}))_{\alpha} \otimes \mathbb{C} t^m$ and
 \[\tilde{\tau}_{\delta_m}=\begin{cases} \mathfrak{h} \otimes \mathbb{C} t^m \oplus I\otimes \mathbb{C}t^m \oplus \mathbb{C}\;\text{ad}\;t^m & \text{if}\; m \in \mathbb{Z}^{n+1}\setminus \text{rad}\;f,\\
 \mathfrak{h}\otimes \mathbb{C}t^m\oplus I \otimes \mathbb{C} t^m\oplus (HC_1 (\mathbb{C}_q))m \oplus (\underset{0\leq i\leq n}\oplus\mathbb{C}t^m\partial_i) &\text{if}\; m \in \text{rad}\;f\setminus \{0\},u \in \mathbb{C}^{n+1},\end{cases}\]
 where $(HC_1(\mathbb{C}_q))_m=\mathbb{C}\langle t^a,t^b\rangle$ with $a+b=m$. \\ 
We define a non-degenerate symmetric bilinear form on $\mathfrak{H}$ by taking the standard form on $\mathfrak{h}$ and extending as $(\mathfrak{h},C_i)=(\mathfrak{h},\partial_i)=(\mathfrak{h},I)=0,\; (C_i,C_j)=(\partial_i,\partial_j)=(C_i,I)=(\partial_i,I)=0$ and $(C_i,\partial_j)=\delta_{ij},\;(I,I)=1$. Similarly we can define a non-degenerate symmetric bilinear form on $\mathfrak{H}^*$. 
For a real root $\gamma$, we define the reflection $r_{\gamma}$ on $\mathfrak{H}^*$ by
\begin{center}
$r_{\gamma}(\lambda)=\lambda-\lambda(\gamma^{\vee})\gamma$, $\lambda\in \mathfrak{H}^*$.
\end{center}
The group generated by $r_{\gamma},\;\gamma\in \bigtriangleup^{re}$ is called the weyl group and is denoted by $W$. Note that the weyl group of type $A$ for the full toroidal Lie algebra is same as this weyl group. It is easy to verify that the form on $\mathfrak{H}$ is $W$-invariant.
\section{Integrable Modules}\label{section 3}
In this section, we introduce integrable modules. Then we conclude that any irreducible integrable module is a highest weight module or a lowest weight module.
\begin{definition}
A module $V$ over $\tilde{\tau}$ is called integrable if
\begin{enumerate}
\item $V=\underset{\lambda \in \mathfrak{H}^*}\oplus V_{\lambda}$, where $V_\lambda=\{v\in V|\;h.v=\lambda (h)v,\;\forall\;h\in \mathfrak{H}\}$.
\item $\text{dim}\; V_\lambda<\infty,$ $\forall \;\lambda\in \mathfrak{h}^*$.
\item For all $\alpha\in \bigtriangleup^{re}$ and $x\in \tilde{\tau}_\alpha$, $x$ acts locally nilpotently on $V$, i.e. for every $v\in V$, there exists $n=n(x,v)\in \mathbb{N}$ such that $x^n.v=0$.
\end{enumerate}
\end{definition}
For an integrable module $V$ of $\tilde{\tau}$, let us denote the set of all weights of $V$ by $P(V)=\{\lambda\in \mathfrak{h}^*\;|\; V_\lambda\neq (0)\}$. For any $\lambda\in P(V)$, $V_\lambda$ is called the weight space of $V$ of weight $\lambda$ and elements of $V_\lambda$ are called weight vectors of weight $\lambda$. The following lemma is standard. 
\begin{lemma}\label{lemma 3.1}
Let $V$ be an irreducible integrable module for $\tilde{\tau}$. Then:
\begin{enumerate}
\item $P(V)$ is $W$-invariant;
\item $dim\;V_\lambda=dim\;V_{w\lambda}$ for all $w\in W$ and $\lambda\in P(V)$.
\item $\lambda(\alpha^{\vee})\in \mathbb{Z}$ for all $\alpha\in \bigtriangleup^{re},\;\lambda\in P(V).$
\item For $\alpha\in \bigtriangleup^{re}$ and $\lambda\in P(V)$ with $\lambda(\alpha^{\vee})>0$, we have $\lambda-\alpha \in P(V)$.
\item $\lambda(C_i)$ is a constant integer for all $\lambda\in P(V)$.
\end{enumerate}
\end{lemma}

Now we take a natural triangular decomposition of $\tilde{\tau}.$
Define \[HC_1^+=\oplus_{\substack {0\leq i\leq n,\\ 0<s_0\in \mathbb{Z},\\ s\in \mathbb{Z}^n}} \mathbb{C}\langle t_i,t_0^{s_0}t^st_i^{-1}\rangle,\]
\[D^+=\{D(u,r),\text{ad}\; t^s|u\in \mathbb{C}^{n+1},r \in \text{rad}\;f,s \notin \text{rad}\;f;r_0,s_0>0\},\]
\newline
\[HC_1^0=\oplus_{\substack{0\leq i\leq n\\ s \in \mathbb{Z}^n}}\mathbb{C}\langle t_i,t^st_i^{-1}\rangle,\]
\[D^0=\{D(u,r),\text{ad}\; t^s|u\in \mathbb{C}^{n+1},r \in \text{rad}\;f,s \notin \text{rad}\;f;r_0,s_0=0\},\]
\newline
\[HC_1^-=\oplus_{\substack {0\leq i\leq n,\\ 0>s_0\in \mathbb{Z},\\ s\in \mathbb{Z}^n}} \mathbb{C}\langle t_i,t_0^{s_0}t^st_i^{-1}\rangle,\]
\[D^-=\{D(u,r),\text{ad}\; t^s|u\in \mathbb{C}^{n+1},r \in \text{rad}\;f,s \notin \text{rad}\;f;r_0,s_0<0\},\]
\newline
\[\tilde{\tau}^+=(\mathfrak{sl}_d(\mathbb{C}))_+\otimes \mathbb{C}_q[t_1^{\pm 1},\dots,t_n^{\pm 1}]\oplus \mathfrak{gl}_d(\mathbb{C})\otimes t_0\mathbb{C}_q[t_0,t_1^{\pm 1},\dots , t_n^{\pm 1}]\oplus HC_1^+ \oplus D^+,\]
\[\tilde{\tau}^-=(\mathfrak{sl}_d(\mathbb{C}))_-\otimes \mathbb{C}_q[t_1^{\pm 1},\dots,t_n^{\pm 1}]\oplus \mathfrak{gl}_d(\mathbb{C})\otimes t_0^{-1}\mathbb{C}_q[t_0^{-1},t_1^{\pm 1},\dots , t_n^{\pm 1}]\oplus HC_1^- \oplus D^-,\]
\[\tilde{\tau}^0=\mathfrak{h}\otimes  \mathbb{C}_q[t_1^{\pm 1},\dots,t_n^{\pm 1}]\oplus I\otimes  \mathbb{C}_q[t_1^{\pm 1},\dots,t_n^{\pm 1}]\oplus HC_1^0\oplus D^0.\]
Let $\bigtriangleup^+=\{\alpha+\delta_m,\delta_{m^{\prime}}|m_0,m_0^{\prime}>0,\text{or}\; m_0=0\; \text{and}\;\alpha >0\}$. Now by Lemma (\ref{lemma 3.1})(5), we assume that
\begin{center}
$\lambda(C_i)=c_i$ for $0\leq i\leq n$, $\forall\; \lambda\in P(V)$ and $c_i\in \mathbb{Z}.$
\end{center}
\begin{theorem}\label{Theorem 3.1}
Let $V$ be an irreducible integrable module for $\tilde{\tau}$ with finite dimensional weight spaces.
\begin{enumerate}
\item If $c_0>0$ and $c_1=c_2=\dots =c_n=0$, then there exists a nonzero element $v\in V$ such that $\tilde{\tau}^+.v=0$.
\item If $c_0<0$ and $c_1=c_2=\dots =c_n=0$, then there exists a nonzero element $v\in V$ such that $\tilde{\tau}^-.v=0$.
\item If $c_0=c_1=\dots =c_n=0$, then there exists a nonzero element $v,w\in V$ such that $(\mathfrak{sl}_d(\mathbb{C}))_+\otimes \mathbb{C}_q.v=0$ and $(\mathfrak{sl}_d(\mathbb{C}))_-\otimes \mathbb{C}_q.w=0$.
\end{enumerate}
\end{theorem}
\begin{proof}
Note that if $\langle t_0^a,t_0^b\rangle\neq 0$ in $V$, then $a+b=0$. Therefore $g_{\text{aff}}=\mathfrak{sl}_d(\mathbb{C})\otimes \mathbb{C}[t_0,t_0^{-1}]\oplus \mathbb{C}C_0\oplus \mathbb{C}\partial_0$ will become a subalgebra of $\tilde{\tau}$. Now the weyl group of $\tilde{\tau}$ is same as the Weyl group of type $A$ for the commuting torus. The proof will follow from Theorem $2.1$ of \cite{[8]}. Note that in order to prove Theorem $2.1$ of \cite{[8]}, they used the properties of weight system of integrable modules and these properties are available in our case also.  
\end{proof}
\begin{lemma}{(Lemma $6.2$(a), \cite{[9]})}
$\mathbb{C}_q\cong \mathbb{C}_q^{\prime}$ iff there exists $A\in GL(n+1,\mathbb{Z})$ such that $q^{\prime}_{ij}=\underset{k,l}\prod q_{kl}^{a_{ki}a_{lj}}$.
\end{lemma}
Let $B=(b_{ij})={(A^t)}^{-1}$. We can extend this isomorphism to an automorphism of $\tilde{\tau}$ as follows (denote this automorphism  by $A$ also)
 
\[A(x(r))=x(A.r),\]
\newline
\[A(\langle t_i,t^mt_i^{-1}\rangle)=\sum_{j=0}^n a_{ji}\langle t_j,t^{(A.m)}t_j^{-1}\rangle,\]   \[0\leq i\leq n,\]
\newline
\[A(t^m\partial_i)=\sum_{j=0}^nb_{ij}t^{(A.m)}\partial_j\]
\newline
\[A(ad\;t^s)=ad\;t^{(A.s)},\]
where $x\in \mathfrak{gl}_d(\mathbb{C}),\;r\in \mathbb{Z}^n,\;m\in \textit{rad}\;f, s\notin \textit{rad}\;f$.\\
Now using this automorphism and Theorem \ref{Theorem 3.1}, we can prove the the following theorem.
\begin{theorem}\label{Theorem 3.2}
Let $V$ be an irreducible integrable module for $\tilde{\tau}$ with finite dimensional weight spaces. Assume $C_i$ acts non-trivially on $V$ for some $0\leq i\leq n$. Then upto automorphism there exists a weight $\lambda \in P(V)$ such that $\lambda+\gamma\notin P(V)$ for all $\gamma\in \bigtriangleup^+$ (or $\lambda-\gamma\notin P(V)$ for all $\gamma\in \bigtriangleup^+$).
\end{theorem}
From now on throughout the paper we will assume $c_0>0$ and $c_1=\dots=c_n=0$. 
By previous theorem $V^+=\{v\in V|\;\tilde{\tau}^+.v=0\}$ is nonzero. 
\begin{lemma}\label{lemma 3.3}
$V^+$ is a $\tilde{\tau}^0$ submodule of $V$. Infact it is irreducible. Further we can say that $V=U(\tilde{\tau}^{-})V^+$. 
\end{lemma}
\begin{proof}
Using PBW theorem and some weight arguments, one can prove this Lemma.
\end{proof}
\begin{lemma}\label{lemma 3.4}
There exists unique $\gamma\in \mathfrak{h}^*$ and $\beta\in \mathbb{C}^n$(not necessarily unique) such that the weights of $V+$ are of the form $\gamma+\delta_{r+\beta}$, where $r\in \mathbb{Z}^n$.
\end{lemma}
\begin{proof}
We can see that $[\mathfrak{h},\tilde{\tau}^0]=0$. Therefore by Lemma \ref{lemma 3.3}, $\mathfrak{h}$ acts scalarly on $V^+$ and hence it will act by a single linear functional. We will denote this functional as $\gamma$. Now we can check that the action of any $r$ degree element of $\tilde{\tau}^0$ on any non zero weight vector of $V^+$ of weight $\mu$ gives us a weight vector of weight $\mu +\delta_r$.
\end{proof}
Since $\partial_i\in \tau^0$ for $1\leq i\leq n$, the highest weight space will be $\mathbb{Z}^n$-graded. So \[V^+=\underset{m\in \mathbb{Z}^n}\oplus V^+(m),\] where $V^+(m)=\{v\in V^+:\partial_i.v=(\lambda(\partial_i)+m_i)v, 1\leq i\leq n\}$.
\section{Action of Central extension part on Highest weight space}\label{section 4}
In this section, we will show that a large part of $HC_1^0$ acts trivially and rest of the part acts associatively on $V^+$. It is easy to see \[HC_1^0=\underset{\substack {s\in \textit{rad}\;f \\ s_0=0}}\oplus\mathbb{C}\;\langle t_0,t^st_0^{-1}\rangle+\underset{\substack {1\leq i\leq n\\s\in \textit{rad}\;f,s_0=0}}\oplus\mathbb{C}\;\langle t_i,t^st_i^{-1}\rangle\]
Now our aim is to prove $HC_1^n=\underset{\substack {1\leq i\leq n\\s\in\text{rad}\;f,s_0=0}}\oplus\mathbb{C}\;\langle t_i,t^st_i^{-1}\rangle$ acts trivially on $V^+$. 
Let $\mathbb{C}_q(n)$ be the Laurent polynomial ring in $n$ non commutating variables $t_1,\dots,t_n$ with the conditions $t_it_j=q_{ij}t_jt_i$. Let $L$ be the Lie algebra $\mathfrak{sl}_d(\mathbb{C})\otimes \mathbb{C}_q(n)\oplus I\otimes [\mathbb{C}_q(n),\mathbb{C}_q(n)]\oplus   HC_1^n\oplus D(n)$, where $D(n)$ be the $\mathbb{C}$-linear span of $\partial_1,\dots \partial _n$. Now we consider a triangular decomposition of $L=L^-\oplus L^0\oplus L^+$, where
$L^+=(\mathfrak{sl}_d(\mathbb{C}))_+\otimes \mathbb{C}_q(n),\;
L^-=(\mathfrak{sl}_d(\mathbb{C}))_-\otimes \mathbb{C}_q(n),$ and 
$L^0=\mathfrak{h}\otimes \mathbb{C}_q(n)\oplus I\otimes [\mathbb{C}_q(n),\mathbb{C}_q(n)]\oplus HC_1^n\oplus D(n).$
 
\begin{theorem}[{\cite{[2]}}, Proposition 2.4]\label{theorem 4.1}
Let $V_1$ be an irreducible integrable module for $L$ with finite dimensional weight spaces. After a suitable co-ordinate change for $\mathbb{C}_q(n)$ we have the followings on $V_1.$
\begin{enumerate}
\item There exists a non-negative integer $k$ and non-zero central operators $z_1,\dots z_k$ of degree $(p_1,0,\dots,0),\dots , (0,\dots, p_k,\dots , 0)$ for some integres $p_1,\dots, p_k$.
\item $k<n$, if $\mathbb{C}_q$ is rational.
\item $\langle t_i,t^rt_i^{-1}\rangle\neq 0$ on $V$ will imply $i\geq k+1$ and $r_{k+1}=\dots =r_{n}=0$.
\end{enumerate} 
\end{theorem}
Note that in \cite{[2]}, the authors expected that if $C_i=\langle t_i, t_i^{-1}\rangle$ acts trivially on any irreducible module of $L$ for $1\leq i\leq n$, then the whole central extension part acts trivially on that module. Here we are giving a proof of it.
\begin{proposition}\label{proposition 4.1}
Let $V_1$ be an irreducible integrable module for $L$ with finite dimensional weight spaces. $k$ be as above theorem. suppose $k\geq 1$ and $C_i=0$ for $i=1,\dots ,n$. Then such a module does not exists.
\end{proposition}
\begin{proof}
The proof will be parallel to Proposition 4.13 of \cite{[4]}. Here we need to construct the Hisenberg Lie algebra differently. Let$\langle t_i,t^mt_i^{-1}\rangle\neq 0$, then by Lemma 2.1, we have $m \in \text{rad}\;f$ and by Theorem \ref{theorem 4.1}, we can see that $i\geq k+1$ and $m_{k+1}=\dots m_n=0$.\\
Since $\mathbb{C}_q(n)$ is rational, choose $N$ be the smallest integer such that $Ne_i\in \text{rad}\;f$. Now consider the Hisenberg algebra $H=\text{span}\;\{ h\otimes t^m t_i^{Nk},h\otimes t_i^{-Nk}, \langle t_i,t^mt_i^{-1}\rangle\;|\;k>0\}$ 
\end{proof}
Let us fix some $i, \;1\leq i\leq n$ and consider the Lie subalgebra $\mathfrak{sl}_d(\mathbb{C})\otimes \mathbb{C}[t_i^{\pm 1}]\oplus \mathbb{C}\partial _i$ of $\tau$. Note that we can consider this loop algebra since we can check $\langle t_i,t_i^mt_i^{-1}\rangle=0$ if $m\neq 0$ and by assumption $\langle t_i,t_i^{-1}\rangle=0$. Let $\theta$ be the highest weight root of $\mathfrak{sl}_d(\mathbb{C})$ and $\theta^{\vee}$ be its corresponding co-root. Since $\mathfrak{h}$ acts scalarly on $V^+$, let $\gamma \in \mathfrak{h}^*$ such that $h.v=\gamma(h)v$ for all $h\in \mathfrak{h}$ and $v\in V^+$. Here $\gamma$ is the restriction of $\lambda$ to $\mathfrak{h}$. Let $W_0,W_i$ be the weyl group of $\mathfrak{sl}_d(\mathbb{C)}$ and loop algebra respectively. Let $r_{\lambda}=min_{h\in \mathbb{Z}(W_0\theta^\vee)}\{\lambda(h):\lambda(h)>0\}\in \mathbb{N}$. Then by 2.4 of \cite{[10]}, we have the following Lemma.
\begin{lemma}\label{lemma 4.1}
For any $r\in \mathbb{Z}$, let us consider $\gamma+r\delta_i$, then there exists $w\in W_i$ such that $w(\gamma+r\delta_i)=\gamma+\bar{r}\delta_i$, where $0\leq \bar{r}<r_{\lambda}$.
\end{lemma}
Let us consider the Lie algebra $\widetilde{\mathfrak{h}}(n)=\mathfrak{h}\otimes \mathbb{C}_q(n)\oplus I\otimes \mathbb{C}_q(n)\oplus HC_1^n$ and $W$ be the Weyl group of the Lie algebra L.
Therefore $W_i \subseteq W$ for $1\leq i\leq n$. Then by Lemma \ref{lemma 4.1}, we have the following corollary
\begin{corollary}
Let $\delta_r=\sum_{i=1}^nr_i\delta_i$, where $r=(r_1,\dots r_n)\in \mathbb{Z}^n$. Let $r_i=\bar{r_i}+s_ir_{\lambda}$, $0\leq \bar{r_i}< r_{\lambda}$. Then there exists $w\in W$ such that $w(\gamma+\delta_r)=\gamma+\delta_{\bar{r}}$, where $\bar{r}=(\bar{r_1},\dots,\bar{r_n})\in \mathbb{Z}^n$. 
\end{corollary}\label{corollary 4.1}
Now we consider $M=\{V^+(l):l\in \mathbb{Z}^n,0\leq l_i<r_{\lambda}, 1\leq i\leq n\}$. Then clearly $M$ is a finite set. Let $N^{\prime}=\textit{Max}\;\{dim\;V^+(l):V^+(l)\in M\}$.

\begin{proposition}\label{proposition 4.2}
The dimensions of the weight spaces of $V^+$ are uniformly bounded by $N^{\prime}$.
\end{proposition}
\begin{proof}
Follows from Corollary \ref{corollary 4.1}.1 .
\end{proof}
Let us take $N=\underset{\substack{l\in \mathbb{Z}^n\\ 0\leq l_i<r_\lambda}}\oplus V^+(l)$. So by Proposition \ref{proposition 4.2}, $N$ is a finite dimensional vector subspace of $V^+$. Now we consider two sets $S_1=\{W|\; W \text{is a}\; L-\text{submodule of}\; V, W\cap V^+\neq 0\}$ and $S_2=\{W|\; W \;\text{is a vector subspace of}\;N\}$. So $S_i$ is partially ordered set with the relation $A\leq B$ if and only if $A\supseteq B$ for $i=1,2$.
\begin{proposition}\label{proposition 4.3}
For every strictly increasing chain of $S_1$, there exists a strictly increasing chain of $S_2$. 
\end{proposition}
\begin{proof}
Let $N_1\supsetneq N_2$ in $S_1$. Therefore there exists $W_i\in S_2$ such that $N_i=U(L)W_i$ for $i=1,2$. Let us assume $W_i^\prime=U(L^0)W_i,\;\widetilde{W_i}=W_i^\prime\cap N,\;\forall\; i=1,2$. Since $N_1\supsetneq N_2$, we can see that $W^\prime_1\supset W^\prime_2$ and hence $\widetilde{W_1}\supset \widetilde{W_2}$. Now we will prove that $U(L^0)W_i=U(L^0)\widetilde{W_i}$ for $i=1,2$.\\
Let $w=\sum u_iw_i\in U(L^0)\widetilde{W_i}$, where $u_i\in U(L^0)$ and $w_i\in \widetilde{W_i}$. Now $\widetilde{W_i}\subset W_i^\prime=U(L^0)W_i$, so $w \in U(L^0)W_i$. Hence we get $U(L^0)\widetilde{W_i}\subset U(L^0)W_i$. Now the converse is trivial because $W_i\subset \widetilde{W_i}$. Now we claim that $\widetilde{W_1}\supsetneq \widetilde{W_2}$. If not, then we have $\widetilde{W_1}=\widetilde{W_2}$. This will imply $U(L)\widetilde{W_1}=U(L)\widetilde{W_2}$ and therefore we get $U(L)W_1=U(L)W_2$ i.e. $N_1=N_2$, a contradiction.
\end{proof}
From Proposition \ref{proposition 4.2} and Proposition \ref{proposition 4.3}, we can say that any increasing chain of $S_1$ has finite length. Let $V_{min}$ be a maximal element of $S_1$. Therefore $V_\text{min}$ is a minimal $L$-submodule of $V$ such that it intersects $V^+$ nontrivially. Note that $V_\text{min}$ is not necessarily an irreducible module for $L$. By the minimality of $V_\text{min}$, we can say that every $L$-submodule of $V_\text{min}$ intersects $V^+$ trivially. So sum of all proper submodule of $V_\text{min}$ is again a proper $L$-submodule of $V_\text{min}$.  So the quotient module of $V_\text{min}$ by this proper submodule is a non-trivial irreducible module for $L$ with finite dimensional weight spaces, where $V_{min}\cap V^+$ goes injectively. Now by Proposition \ref{proposition 4.1}, we can say that on this module $HC_1^n$ acts trivially and hence it acts trivially on $V_\text{min}\cap V^+$.
\begin{proposition}\label{proposition 4.4}
$HC_1^n$ acts trivially on $V^+$.
\end{proposition}
\begin{proof}
By previous discussion $M^{\prime}=\{v\in V^+:HC_1^n.v=0\}$ is a nonzero subspace of $V^+$. Using Lemma \ref{lemma 2.1}, we can see that $[D(u,r),\langle t_i,t^mt_i^{-1}\rangle]=(u,m)\sigma(r,m)\langle t_i,t^{r+m}t_i^{-1}\rangle+u_i\sigma(r,m)\sum_{j=0}^nr_j\langle t_j,t^{m+r}t_j^{-1}\rangle$. Therefore $[D^0,HC_1^n]\subset HC_1^n$. We know $HC_1^n$ commutes with $\mathfrak{sl}_d(\mathbb{C})\otimes \mathbb{C}_q(n)\oplus I\otimes \mathbb{C}_q(n)$. Hence $M^{\prime}$ is a nonzero $\tau^0$ sub-module of $V^+$. Now by irreducibility of $V^+$, we have the Lemma.
\end{proof}
Now we will prove that $\langle t_0,t^mt_0^{-1}\rangle$ acts associatively on the highest weight space for $m\in \textit{rad}\;f, m_0=0$. Let us recall that $\lambda(C_0)=c_0>0$. This is a crucial thing which will be used in order to prove the associative action.
\begin{lemma}\label{lemma 4.2}
Let $\langle t_0,t^mt_0^{-1}\rangle.v=0$ for some $v\in V^+$ and $m\in \textit{rad}\; f,m_0=0$. Then $\langle t_0,t^mt_0^{-1}\rangle$ will act locally nilpotently on $V^+$.
\end{lemma}
\begin{proof}
Consider the set $M^{\prime}=\{w\in V^+:\textit{there exists}\; k\in \mathbb{N}\;\text{such that}\;(\langle t_0,t^mt_0^{-1}\rangle)^k.w=0\}$. By assumption this set is nonzero. We know that \[[\mathfrak{h}\otimes \mathbb{C}_q(n)\oplus I\otimes \mathbb{C}_q(n),\langle t_0,t^mt_0^{-1}\rangle]=0.\]
Now using Proposition \ref{proposition 4.4}, we can see that \[[D(u,r),\langle t_0,t^mt_0^{-1}\rangle]=(u,m)\sigma(r,m)\langle t_0,t^{m+r}t_0^{-1}\rangle.\]
Now using this Lie bracket, we will have the following equation
 
\[(\langle t_0,t^mt_0^{-1}\rangle)^{k+1}.D(u,r)=D(u,r).(\langle t_0,t^mt_0^{-1}\rangle)^{k+1}\]\[-(k+1)(u,m)\sigma(r,m)\langle t_0,t^{m+r}t_0^{-1}\rangle(\langle t_0,t^mt_0^{-1}\rangle)^k.\]
Using this equation we can see that $D(u,r)M^{\prime}\subseteq M^{\prime}$ for $r\in \textit{rad}\;f,r_0=0$. Now we have $[t^r\partial_0,\langle t_0,t^mt_0^{-1}\rangle ]=0$ for $r\in \textit{rad}\;f,r_0=0$ and $[ad\;t^s,\langle t_0,t^mt_0^{-1}\rangle]=0$ for $s\notin \textit{rad}\;,s_0=0$. Therefore $M^{\prime}$ is a nonzero $\tau^0$-submodule of $V^+$ and hence by irreducibility of $V^+$, we have the Lemma. 
\end{proof}
As in Section 5 of \cite{[12]}, we can prove the following Lemmas.
\begin{lemma}
Let $\langle t_0,t^mt_0^{-1}\rangle$ acts locally nilpotently on $V^+$ for all $m\in \textit{rad}\;f\setminus\{0\}, m_0=0$. Then $\langle t_0,t^mt_0^{-1}\rangle$ acts trivially on $V^+$ for all $m\in \textit{rad}\;f\setminus\{0\}, m_0=0$.
\end{lemma}
\begin{lemma}\label{lemma 4.4}
$\langle t_0,t^mt_0^{-1}\rangle$ acts injectively or trivially on $V^+$, for $m_0=0,m\in \textit{rad}\;f$.
\end{lemma}
Now by our assumption $c_0\neq 0$. If there exists  $r\in \textit{rad}\; f\setminus\{0\}$ with $r_0=0$ such that $\langle t_0,t^rt_0^{-1}\rangle$ is not injective, then by Lemma \ref{lemma 4.4}, we can see that $\langle t_0,t^mt_0^{-1}\rangle$ acts trivially for all $m\in \textit{rad}\;f\setminus\{0\}$ with $m_0=0$. Let $\theta$ be the highest root of $\mathfrak{sl}_d(\mathbb{C})$ and $x_{\theta}$ and $x_{-\theta}$ be elements of $(\mathfrak{sl}_d(\mathbb{C}))_{\theta}$ and $(\mathfrak{sl}_d(\mathbb{C}))_{-\theta}$  respectively such that $tr(x_\theta,x_{-\theta})=1$. Let $n_0$ be the smallest positive integers such that $t_0^{n_0}\in Z(\mathbb{C}_q)$. Now it is easy to see that $x_{\theta}\otimes t_0^{n_0}$ and $x_{-\theta}\otimes t_0^{-n_0}$ and $h_{\theta}+n_0\;C_0$ form a $\mathfrak{sl}_2$ copy. Assume $n_1$ be the order of $q_{1,0}=q_{0,1}^{-1}$. Now let us consider the loop of this algebra $\mathfrak{sl}_2\otimes \mathbb{C}[t_1^{\pm n_1}]$ with the bracket:
\[[(x_{\theta}\otimes t_0^{n_0})\otimes t_1^{r_1},(x_{-\theta}\otimes t_0^{-n_0})\otimes t_1^{s_1}]=h_{\theta}\otimes t_1^{r_1+s_1}+n_0\;C_0\otimes t_1^{r_1+s_1},\]
where $r_1,s_1\in n_1\mathbb{Z}.$
In our setting the Lie bracket is becoming $[x_{\theta}(n_0e_0+r_1e_1),x_{-\theta}(-n_0e_0+s_1e_1)]=h_{\theta}(r_1e_1+s_1e_1)+n_0\langle t_0,t_1^{r_1+s_1}t_0^{-1}\rangle+r_1\langle t_1,t_1^{r_1+s_1}t_1^{-1}\rangle.$ Since $r_1\langle t_1,t_1^{r_1+s_1}t_1^{-1}\rangle$ acting trivially on $V^+$, we can identify $n_0\;C_0\otimes t_1^{r_1+s_1}$ with $n_0\langle t_0,t_1^{r_1+s_1}t_0^{-1}\rangle$. Now since $V$ is an integrable module for $\tilde{\tau}$ in particular for $L(\mathfrak{sl}_2)=\mathfrak{sl}_2\otimes \mathbb{C}[t_1^{\pm n_1}]$ and every element of $V^+$ is a highest weight vector for $L(\mathfrak{sl}_2)$. Now using \cite{[13]} as in Proposition 3.10 of \cite{[11]}, we can say that $\langle t_0,t^mt_0^{-1}\rangle.v\neq 0$ for all $m\in \textit{rad}\;f,m_0=0$ and $v\in V^+$.
\vspace{5mm}

Now let us concentrate at the restricted $n\times n$ matrix $q^{\prime}$ of $q$ obtained by omitting the first row and first column of $q$. Let $\sigma^{\prime},f^{\prime}$ be the restriction function of $\sigma,f$ on $\mathbb{Z}^n\times \mathbb{Z}^n$. Therefore $\mathbb{C}_{q^{\prime}}$ is the Laurent polynomial ring in $n$ non-commutating variables $t_1,t_2,\dots t_n$ with the matrix $q^{\prime}$.
\vspace{5mm}

Using Theorem 4.5 of \cite{[14]}, as in section 2 of \cite{[15]}, upto isomorphism of $\mathbb{C}_{q^\prime}$, we may assume that $q^\prime$ is in simple form. Therefore as in \cite{[15]}, $\textit{rad}\;f^\prime=\xi_1\mathbb{Z}\oplus\dots \oplus \xi_n\mathbb{Z}$.  One of the main advantage of this simple form we are getting is $\sigma\prime(a,b)=\sigma^\prime(b,a)=1$ for all $a\in \textit{rad}\;f^\prime$ and $b\in \mathbb{Z}^n$. Since $V^+$ is $\mathbb{Z}^n$-graded $\tau^0$-module, we can write $V^+=\underset{r\in \textit{rad}\;f^\prime}\oplus V^+_r$, where $V^+_r=\underset{r_i\leq k_i<\xi_i+r_i}\oplus V^+(k)$ for $k\in \mathbb{Z}^n$, therefore it will be $\textit{rad}\;f^\prime$-graded also. 
Now the following proposition can be proved as in Proposition 3.10 of \cite{[11]}.

\begin{proposition}\label{proposition 4.5}
\begin{enumerate}
\item $\langle t_0,t^rt_0^{-1}\rangle \langle t_0,t^st_0^{-1}\rangle=c_0\;\langle t_0,t^{r+s}t_0^{-1}\rangle$ on $V^+$ for all $r,s\in \textit{rad}\; f^\prime$.
\item $dim\; V^+(k)=dim\; V^+(k+r)=p_k$ for all $r\in \textit{rad}\;f^\prime$. Suppose $v_1(k),\dots,v_{p_k}(k)$ is a basis for $V^+(k)$, where $0\leq  k_i<\xi_i$. Let $v_i(k+r)=\frac{\langle t_0.t^rt_0^{-1}\rangle}{c_0}v_i(k)$,  $\forall\; i,\forall\; r\neq 0$. Then $v_1(k+r),\dots , v_{p_k}(k+r)$ is a basis for $V^+(k+r)$.
\item For $0\neq r\in \textit{rad}\;f^\prime,$ $\langle t_0,t^st_0^{-1}\rangle(v_1(k+r),\dots v_{p_k}(k+r))=c_0(v_1(k+r+s),\dots v_{p_k}(k+r+s))$.
\item $0\neq r\in \textit{rad}\;f^\prime$, we have $t^s\partial_0(v_1(k+r),\dots v_{p_k}(k+r))=\lambda (\partial_0)(v_1(k+r+s),\dots v_{p_k}(k+r+s))$.
\end{enumerate}
\end{proposition}
Note that we can not prove the associative action of $h(a)$ for all $a\in \mathbb{Z}^n$, since  Lemma \ref{lemma 4.2} is not true for $h(a)$, when $a\notin \textit{rad}\;f^\prime$.
\section{Complete reducibility}\label{section 5}
In this section we will prove that a quotient of the highest weight space is completely reducible for a subalgebra of $\tau^0$.\\
Let us assume $D(u,r)=\underset{1\leq i\leq n}\sum u_it^r\partial_i\in D^0$. Let $I(u,r)=D(u,r)-D(u,0)$ for $u\in \mathbb{C}^n, r\in \textit{rad}\;f^\prime$. Consider the subspace $\textit{span}\;\{I(u,r),ad\;t^s:u\in \mathbb{C}^n,r\in \textit{rad}\;f^\prime,s\in \mathbb{Z}^n\setminus \textit{rad}f^\prime\}$. We can see that this will be a Lie subalgebra of $Der(\mathbb{C}_{q^\prime})$ with the Lie brackets
\[[I(u,r),I(v,s)]=I(w,r+s)+(v,r)I(u,r)-(u,s)I(v,s),\]
\[[I(u,r),ad\;t^a]=(u,a)(ad\;t^{r+a}-ad\;t^a),\]
where $u,v\in \mathbb{C}^n,\;r,s\in \textit{rad}\;f^\prime, a\notin \textit{rad}\;f^\prime$.
\vspace{5mm}

Therefore $(\mathfrak{h}\otimes \mathbb{C}_{q^\prime}\oplus I\otimes\mathbb{C}_{q^\prime}\oplus HC_1^0\underset{r\in \textit{rad}\;f^\prime}\oplus \mathbb{C}t^r\partial_0)\rtimes \textit{span}\;\{I(u,r),ad\;t^s:u\in \mathbb{C}^n,r\in \textit{rad}\;f^\prime,s\in \mathbb{Z}^n\setminus \textit{rad}f^\prime\} $ will become a Lie subalgebra of $\tilde{\tau}^0$.
Consider the subspace $W=\textit{span}\{\langle t_0,t^mt_0^{-1}\rangle.v-v:v\in V^+,m\in \textit{rad}\;f^\prime\}$. Then one can easily check that this will be a submodule of $V^+$ for the above Lie algebra. Now let us take the quotient $\tilde{V}^+=V^+/W$, this will be a finite dimensional module for the above Lie algebra. It is easy to see from Proposition \ref{proposition 4.5}, that $\langle t_0,t^mt_0^{-1}\rangle$ and $t^r\partial_0$ acts scalarly on $\tilde{V}^+$, where $m,r \in \text{rad}\;f^\prime$. Therefore we can identify the above Lie algebra with $\mathfrak{L}=(\mathfrak{h}\otimes \mathbb{C}_{q^\prime}\oplus I\otimes\mathbb{C}_{q^\prime})\rtimes \textit{span}\;\{I(u,r),ad\;t^s:u\in \mathbb{C}^n,r\in \textit{rad}\;f^\prime,s\in \mathbb{Z}^n\setminus \textit{rad}f^\prime\}$ and $\tilde{V}^+$ is a module for $\mathfrak{L}$. Infact in this section we will show that $\tilde{V}^+$ is completely reducible module for $\mathfrak{L}$. 
\vspace{5mm}

Let $\lambda\in P(V)$ be as of Theorem \ref{Theorem 3.2} and $\alpha_i=\lambda(\partial_i)$ for $1\leq i\leq n$. Let $\alpha=\sum \alpha_ie_i\in \mathbb{C}^n$ and assume $V_1$ be any $\mathfrak{L}$-module. Then we can define a $\tilde{\tau}^0$-module structure on $L(V_1)=V_1\otimes \mathbb{C}_q$.
\[h(a).v\otimes t^b=(h(a)v)\otimes t^{a+b},\]
\[I(a).v\otimes t^b=(I(a)v)\otimes t^{a+b},\]
\[D(u,r).v\otimes t^b=(I(u,r)v)\otimes t^{b+r}+(u,b+\alpha)v\otimes t^{b+r},\]
\[\langle t_0,t^mt_0^{-1}\rangle.v\otimes t^a=c_0v\otimes t^{m+a},\]
\[t^r\partial_0.v\otimes t^a=\lambda (\partial_0)v\otimes t^{a+r},\]
\[ad\;t^s.v\otimes t^a=(ad\;t^sv)\otimes t^{a+s},\]
where $v\in V_1,\; a,b\in \mathbb{Z}^n,\;m,r \in \textit{rad}\;f{\prime},\; s\notin \textit{rad}\;f^\prime$.
\vspace{5mm}

By the above discussion $L(\tilde{V}^+)$ has a $\tilde{\tau}^0$-module structure. We will show that $V^+$ is contained in $L(\tilde{V}^+)$ as an $\tilde{\tau}^0$-module. For $v_a\in V^+(a)$, let $\bar{v}_a$ be the image in $\tilde{V}^+=V^+/W$. Let $\tilde{\phi}:V^+\rightarrow L(\tilde{V}^+)$ be $\tilde{\phi}\;(v_a)=v_a\otimes t^a$ for $a\in \mathbb{Z}^n$.
\begin{lemma}
$\tilde{\phi}$ is an $\tilde{\tau}^0$-module map.
\end{lemma}
\begin{proof}
Consider $\tilde{\phi}(ad\;t^s.v_a)=(ad\;t^s.v_a)\otimes t^{a+s}$. Now let us see the action of $ad\;t^s$ on $v_a\otimes t^a$. $ad\;t^s.(v_a\otimes t^a)=(ad\;t^s.v_a)\otimes t^{a+s}$. Rest of the calculation will be same as Lemma 8.2 of \cite{[7]}.
\end{proof}
Now one can proceed same as Section 8 of \cite{[7]} and prove the following theorem.
\begin{theorem}\label{theorem 5.1}

\begin{enumerate}
\item $V^+\cong L(\widetilde{V}^+)(\bar{0})$ as $\tau^0$-module. 
\item $\widetilde{V}^+$ is $\mathbb{Z}^n/\textit{rad}\;f^\prime$-graded irreducible module for $\mathfrak{L}$. 

\item $\tilde{V}^+$ is completely reducible as $\mathfrak{L}$-module. All its irreducible components are isomorphis as $\mathfrak{h}\otimes Z(\mathbb{C}_q^\prime)\oplus I\otimes Z(\mathbb{C}_q^\prime))\rtimes \textit{span}\;\{I(u,r)|u\in \mathbb{C}^n, r\in \textit{rad}\;f^\prime\}$.
\end{enumerate}
 
\end{theorem} 
\section{identification of highest weight space}\label{section 6}
Now let us look at realization $\mathbb{C}_{q^\prime}$ as loop algebra. As in \cite{[15]}, let $q^\prime$ is in simple form, so we can assume that $q_{2i, 2i-1}=q_i,\; q_{2i-1,2i}=q_i^{-1}$, for $1\leq i \leq z$, and other entries of q are all $1$, where $z\in \mathbb{N}$ with $2z\leq n$ and with the orders $k_i$ of $q_i$, $1\leq i \leq z$ as roots of unity satisfy $k_{i+1}|k_i,\; 1\leq i< z$.  \\
For $1\leq l \leq n$, let 
$\xi_{l}=\begin{cases} k_ie_{2i-1}, & \textit{if}\; l=2i-1\leq 2z,\\ k_ie_{2i} & \textit{if}\; l=2i \leq 2z, \\ e_l, & \textit{if}\; l>2z.
\end{cases}$\\
Then $\{\xi_1, \dots ,\xi_n\}$ is a $\mathbb{Z}$-basis of the subgroup $\textit{rad} \;(f^\prime)$.  Denote $\Gamma=\mathbb{Z}^n/\textit{rad}\;f^\prime$ and $\Gamma_0=\{r\in \mathbb{Z}^n|0<n_{2i-1}\leq k_i,\;0<n_{2i}\leq k_i, 1\leq i \leq z,\textit{and}\; n_l=0,\;2z<l\leq n\}$. we can see that $\Gamma_0$ generates $\Gamma$.\\
Now  $\mathfrak{I}=\textit{span}\;\{t^{a+r}-t^a|a\in \mathbb{Z}^n,\; r\in \textit{rad}\;f^\prime\}$ will become ideal of the associative algebra $\mathbb{C}_{q^\prime}$ and we have $\mathbb{C}_{q^\prime}/\mathfrak{I}\cong \otimes _{i=1}^zM_{k_i}(\mathbb{C})\cong M_N(\mathbb{C})$ as associative algebra, where $N=\prod_{i=1}^zk_i$. Therefore as in \cite{[15]}, the general linear Lie algebra $\mathfrak{gl}_N=M_N(\mathbb{C})$ will become a $\Gamma$-graded Lie algebra with $\mathfrak{gl}_N=\bigoplus_ {a\in \Gamma}\mathbb{C} X^{\bar{a}}$, where $X$ is defined in \cite{[15]}. From the construction of $X$, we can easily see that $X^r$ is the identity matrix of $\mathfrak{gl}_N$ for any $r\in \textit{rad}\;f^\prime$. Here the Lie brackets are given by $[X^a, X^b]=(\sigma(a,b)-\sigma(b,a))X^{a+b}$. Infact we can see that $\mathfrak{sl}_N=\underset{\bar{a}\in \Gamma\setminus{\bar{0}}}\oplus\mathbb{C}X^{\bar{a}}$.

\vspace{5mm}

Let $A=\mathbb{C}[x_1^{\pm 1},\dots x_n^{\pm 1}]$ be the commuting Laurent polynomial ring in $n$ variables $x_i$. Now consider the multi-loop algebra $\mathfrak{gl}_N\otimes A$ with the Lie bracket $[X^{\bar{a}}\otimes x^r,X^{\bar{b}}\otimes x^s]=(\sigma(a,b)-\sigma(b,a))X^{\bar{a}+\bar{b}}\otimes x^{r+s}$ for $\bar{a},\bar{b} \in \Gamma_0$ and $r,s\in \mathbb{Z}^n$. Now we recall a result of \cite{[15]}, from which we can realize all finite dimensional simple modules for $\mathbb{C}_{q^\prime}$.
\begin{lemma}{(Lemma 2.3,\cite{[15]})}
As associative algebras, \[\mathbb{C}_{q^\prime}\cong  \underset{a\in \mathbb{Z}^n}\bigoplus (\mathbb{C}X^a\otimes x^a),\]  
where the right hand side is a $Z^n-$graded subalgebra of $\mathfrak{gl}_N\otimes A$.
\end{lemma}
Recall that $\mathfrak{gl}_N$ is $\Gamma$-graded with each graded component is one dimensional. For $\bar{a}\in \Gamma$, we have $(\mathfrak{gl}_N)_{\bar{a}}=\mathbb{C}X^{\bar{a}}$. Here $X^{a}=\otimes_{i=1}^zX_{2i-1}^{\bar{a}_{2i-1}}X_{2i}^{\bar{a}_{2i}}$, where $X_{2i-1}=\sum_{j=1}^{k_i}q_i^{j-1}E_{j,j}$ and $X_{2i}=\sum_{j=1}^{k_i-1}E_{j,j+1}+E_{k_i,1}$ are generators for the associative algebra $M_{k_i}(\mathbb{C})$, $q_i$ is primitive $k_i$-th root of unity. Let us define $\sigma_{2i-1},\sigma_{2i}:M_{k_i}(\mathbb{C})\rightarrow M_{k_i}(\mathbb{C})$ as $\sigma_{2i-1}(X_{2i-1})=q_iX_{2i-1},\;\sigma_{2i-1}(X_{2i})=X_{2i}$ and $\sigma_{2i}(X_{2i-1})=X_{2i-1},\;\sigma_{2i}(X_{2i})=q_iX_{2i}$. One can easily check that $\sigma_{2i-1},\;\sigma_{2i}$ are associative algebra automorphisms of $M_{k_i}(\mathbb{C})$ of order $k_i$ and hence they will become Lie algebra automorphisms of $\mathfrak{gl}_N(\mathbb{C})$. Now the common eigen spaces of $\sigma_i$ for $1\leq i\leq 2z$ coincide with the graded components of $\mathfrak{gl}_N(\mathbb{C})$.  Also we can see that the group generated by $\{\sigma_i|1\leq i\leq 2z\}$ is same as $\Gamma$. So if we take the  $\Gamma$-actoin on $A$ by 
\[\sigma_i.f(x_1,\dots,x_n)=f(t_1,\dots,t_{i-1},q_i^{-1}t_i,\dots ,t^n),\; f\in A.\] 
We can see with the notations of \cite{[16]}, $M((\mathbb{C}^*)^n,\mathfrak{gl}_N)^{\Gamma}=\underset{a\in \mathbb{Z}^n}\bigoplus (\mathbb{C}X^a\otimes x^a)$.
Now let us look at finite dimensional simple modules of $\mathbb{C}_{q^\prime}$. Let $W$ be a finite dimensional simple module for $\mathbb{C}_{q^\prime}$. Therefore by above isomorphism it will be an irreducible module for $\underset{a\in \mathbb{Z}^n}\bigoplus (\mathbb{C}X^a\otimes x^a)$. 
\vspace{5mm}

For any  $a=(a_1,\dots, a_n)\in (\mathbb{C}^*)^n$, we define $\xi (a)=(a_1^{\xi_1},\dots, a_n^{\xi_n})$. 
\begin{lemma}{(Corollary 5.4, \cite{[16]})}\label{lemma 6.2}
Let $W$ be a finite dimensional simple module for the twisted loop algebra $\underset{a\in \mathbb{Z}^n}\bigoplus (\mathbb{C}X^a\otimes x^a)$. Then there exists $b_1,\dots ,b_r\in (\mathbb{C}^*)^n$ and $W_{(\lambda_1,a_1)},\dots ,W_{(\lambda_r,a_r)}$ irreducible modules for $\mathfrak{gl}_N$ such that $W\cong W_{(\lambda_1,a_1)}(b_1)\otimes \dots \otimes W_{(\lambda_r,a_r)}(b_r)$, where $\xi (b_i)\neq \xi(b_j)$ whenever $i\neq j$.
\end{lemma}
If $W$ is a finite dimensional simple module for $[\mathbb{C}_{q^{\prime}},\mathbb{C}_{q^{\prime}}]$, then we can make it a simple module for $\mathbb{C}_{q^{\prime}}$ by giving trivial action of $Z(\mathbb{C}_{q^{\prime}})=\underset{r\in \text{rad}\;f^{\prime}}\oplus I_N\otimes x^r$ on the modules. In that case $W$ will be same as Lemma \ref{lemma 6.2}, with $a_i=0$ for $1\leq i\leq r$ with $r\in \mathbb{N}$.
Now consider the subspace $(\mathfrak{h}\otimes \mathbb{C}_{q^\prime}\oplus I\otimes \mathbb{C}_{q^\prime})\rtimes \textit{span}\{ ad\;t^s|s \notin \textit{rad\;}(f^\prime)\}$ of $\mathfrak{L}$. We can easily check that this is an ideal of $\mathfrak{L}$. We would like to give a Lie algebra isomorphism between the above ideal and direct sum of twisted loop algebras.  Now let us consider the Lie algebra $L=\underset{1\leq i\leq d}\oplus L_i\oplus (\underset{a\in \mathbb{Z}^n\setminus{\text{rad}\;f^{\prime}}}\oplus(\mathfrak{sl}_N)_{\bar{a}}\otimes x^a)$, where $L_i=\underset{a\in \mathbb{Z}^n}\oplus(\mathfrak{gl}_N)_{\bar{a}}\otimes x^a=\oplus_{a\in \mathbb{Z}^n}\mathbb{C}(X_i^{\bar{a}}\otimes x^a)$ and $\mathfrak{sl}_N=\oplus_{a\in \Gamma\setminus{\bar{0}}}\mathbb{C}Y^{\bar{a}}$. Here $X_i=Y=X$ for $1\leq i\leq d$.

\begin{proposition}\label{proposition 6.1}
$L\cong (\mathfrak{h}\otimes \mathbb{C}_{q^\prime}\oplus I\otimes \mathbb{C}_{q^\prime})\rtimes \textit{span}\{ ad\;t^s|s \notin \textit{rad\;}(f^\prime)\}$

\end{proposition}
\begin{proof}
We can see that as vector space $\mathfrak{h}\otimes \mathbb{C}_{q^\prime}\oplus I\otimes \mathbb{C}_{q^\prime}\rtimes \textit{span}\{ ad\;t^s|\;s \notin \textit{rad\;}(f^\prime)\}=\underset{1\leq i\leq d}\oplus (E_{ii}\otimes \mathbb{C}_{q^\prime})\oplus \textit{span}\{ ad\;t^s-I\otimes t^s|\;s \notin \textit{rad\;}(f^\prime)\}$. One can easily prove that each component of the direct sum commute with every other component. Now the mapping $X_i^{\bar{a}}\otimes x^a\mapsto E_{ii}\otimes t^{a},\; Y^{\bar{a}}\otimes x^a\mapsto ad\;t^{a}-I\otimes t^{a}$ will give the isomorphism.  
\end{proof}

Let $W$ be a finite dimensional simple module for $L$. Then using Lemma 3.13 of  \cite{[17]}, we can see that $W\cong (\underset{1\leq i \leq d}\otimes W^i)\otimes W^{d+1}$, where each $W^i$ is irreducible $\underset{a\in \mathbb{Z}^n}\oplus(\mathfrak{gl}_N)_{\bar{a}}\otimes x^a $-module for $1\leq i \leq d$ and $W^{d+1}$ is irreducible module for $(\underset{a\in \mathbb{Z}^n\setminus{\text{rad}\;f^{\prime}}}\oplus(\mathfrak{sl}_N)_{\bar{a}}\otimes x^a) $.  Therefore every finite dimensional simple module of $L$ will be tensor product of $(d+1)$ simple modules for $\mathbb{C}_{q^{\prime}}$, where  on the $(d+1)$-th module $Z(\mathbb{C}_{q^{\prime}})$ acts trivially.

\vspace{5mm}

Since we have Theorem \ref{theorem 5.1}, in order to identify the highest weight space as $\tilde{\tau^0}$-module, it is essential to realize finite dimensional simple modules for $\mathfrak{L}$. Let $(\tilde{V}, \pi)$ be such module. Now by Proposition 19.1(b) of \cite{[18]}, we can say that $\pi (\mathfrak{L})$ is reductive with atmost one dimensional center. Therefore $\mathfrak{g}_1=\pi(L)$ being an ideal of $\pi(\mathfrak{L})$, there exists a unique complement $\mathfrak{g}_2$ of $\mathfrak{g}_1$ in $\pi(L)$. So again using \cite{[17]}, we can say that $\tilde{V}\cong V_1\otimes V_2$, where $V_i$ is irreducible module for $\mathfrak{g}_i$, $i=1,2$. Hence from the previous discussion we can see $V_1$ is tensor product of $(d+1)$ irreducible modules for $\mathbb{C}_{q^{\prime}}$, where  on the $(d+1)$-th module $Z(\mathbb{C}_{q^{\prime}})$ acts trivially. So from Lemma \ref{lemma 6.2} and previous discussion we say that $V_1$ is a module for $\oplus\; \mathfrak{gl}_N$ (finitely many copies). Let $\pi_i$ be the $i$-th projection map of $\pi$ on the $i$-th piece of $\oplus \mathfrak{gl}_N$. We claim that $\pi_i(I(u,r))=0.$ Consider the Lie bracket $[I(u,r), X_j^{\bar{a}}\otimes x^{a+s}]=(u,a+s)(X_j^{\bar{a}}\otimes x^{r+a+s}-X_j^{\bar{a}}\otimes x^{a+s})$ for $u\in \mathbb{C}^n,\;a\in \mathbb{Z}^n$ and $r,s\in\text{rad}\;f^{\prime}$. Applying $\pi_i$ both side we will get $[\pi_i(I(u,r)),b_i^{a+s}X^{\bar{a}}]=(u,a+s)b_i^{a+s}(b_i^{r}-1)X_j^{\bar{a}}$ for the $i$-th evolution point $b_i\in (\mathbb{C}^*)^n$. Now canceling $b_i^{a+s}$ from both side we can see that the Lie bracket of LHS is independent of $s$. Similar things hold for the Lie brackets $[\pi_i(I(u,r)), Y^{\bar{a}}\otimes x^{a+s}]$ for $a\in \mathbb{Z}^n\setminus{\text{rad}\;f^{\prime}},\;r,s\in\text{rad}\;f^{\prime} $. So $\pi_i(I(u,r))=0$ for all projection, hence we have $\pi(I(u,r))\subset \mathfrak{g_2}$. Then $V_2$ is irreducible module for $\textit{span}\{I(u,r)|\;u\in \mathbb{C}^n, r\in \textit{rad}\;f^\prime\}$. From Proposition 3.1 and Proposition 3.3 of  \cite{[19]}, we can say that $V_2$ is irreducible module for $\mathfrak{gl}_n$. 
\vspace{5mm}

From Theorem \ref{theorem 5.1}, let us assume $\widetilde{V}^+=\oplus_{i=1}^kM_i$, where $M_i$ is irreducible module for $\mathfrak{L}$, $1\leq i\leq k$ and $k\in \mathbb{N}$. Therefore from previous discussion $M_i\cong M_i^1\otimes M_i^2$, where $M_i^1$ is irreducible module for $\oplus \mathfrak{gl}_N$ and $M_i^2$ irreducible module for $\mathfrak{gl}_n$. Again from Theorem \ref{theorem 5.1}, we can see that $M_i^2\cong M_j^{2}(=M_2, \textit{say})$ as $\textit{span}\; \{I(u,r)|u\in \mathbb{C}^n,r\in \textit{rad}\;f^\prime\}$ module, for $1\leq i,j \leq k$. So we can take $\widetilde{V}^+\cong M_1\otimes M_2$, where $M_1=\sum_{i=1}^{k}M_i^1$ is $\Gamma$-graded irreducible module for $\underset{d+1-\textit{copies}}\bigoplus \mathbb{C}_{q^\prime}$ with trivial action of $Z( \mathbb{C}_{q^\prime})$ from the last copy.
\vspace{7mm}

Let us recall $\tilde{\tau}^0=\mathfrak{h}\otimes \mathbb{C}_{q^\prime}\oplus I\otimes \mathbb{C}_{q^\prime}\bigoplus \underset{m\in \textit{rad}\;f^\prime}\oplus \langle t_0,t^mt_0^{-1}\rangle \bigoplus \underset{m\in \textit{rad}\;f^\prime}\oplus t^m\partial _0\oplus \textit{Der}(\mathbb{C}_{q^\prime})$, where $D^0$ is identified with $\underset{m\in \textit{rad}\;f^\prime}\oplus t^m\partial _0\oplus \textit{Der}(\mathbb{C}_{q^\prime})$. Now we assume $V_1$ be a $\Gamma$-graded finite dimensional irreducible module for $L$ and $V_2$ is finite dimensional irreducible for $\mathfrak{gl}_n$. Now we are going to define a $\tilde{\tau}^0$ action on $L(V_1\otimes V_2)=V_1\otimes V_2\otimes \mathbb{C}_{q^\prime}$.
\[x\otimes t^k.(v_1\otimes v_2\otimes t^l)=(x\otimes t^k.v_1)\otimes v_2\otimes t^{k+l},\]
\[D(u,r).(v_1\otimes v_2\otimes t^l)=v_1\otimes (\sum_{i,j}u_ir_jE_{ji}v_2)\otimes t^{l+r}+(u,k+\alpha)v_1\otimes v_2\otimes t^{l+r},\]
\[<t_0,t^rt_0^{-1}>.(v_1\otimes v_2\otimes t^l)=c_0(v_1\otimes v_2\otimes t^{l+r}),\]
\[t^r\partial_0.(v_1\otimes v_2\otimes t^l)=\lambda(\partial_0)(v_1\otimes v_2\otimes t^{l+r}),\]
\[ad\;t^s.(v_1\otimes v_2\otimes t^k)=(ad\;t^s.v_1)\otimes v_2\otimes t^{k+s},\]
where $x\otimes t^k\in L$ , $k,l\in \mathbb{Z}^n$ and $r\in \textit{rad}\;f^\prime$.\\
One can check that with the above action $L(V_1\otimes V_2)$ is $\tilde{\tau}^0$ module. $V_1=\underset{\bar{k}\in \Gamma}\bigoplus V_{1,\bar{k}}$ is $\Gamma$-graded which is compatible with the $\Gamma$-gradation of $L$. Therefore we can see that $M=\underset{k\in \mathbb{Z}^n}\bigoplus V_{1,\bar{k}}\otimes V_2\otimes t^k$ is an irreducible module for $\tilde{\tau}^0$. Consider the induced module $U(\tilde{\tau})M$, where $\tilde{\tau}^+$ acts trivially. Then it has a maximal submodule (Say $M^{\textit{Rad}}$) which meets $M$ trivially. Then $U(\tilde{\tau})M/M^{\textit{Rad}}$ is irreducible module for $\tilde{\tau}$.
\begin{theorem}\label{theorem 6.1}
Let $V$ be an irreducible integrable $\tilde{\tau}$ module with finite dimensional weight spaces, with $<t_0,t^{-1}_0>$ acting as $c_0$ and $<t_i,t^{-1}_i>$ acts trivially for $1\leq i\leq n$. Then $V\cong U(\tilde{\tau})M/M^{\textit{Rad}}$.
\end{theorem}

\begin{remark}
Using Proposition \ref{proposition 4.1} and Proposition \ref{proposition 6.1}, one can settle the level zero modules for $\tau$ added with zero degree derivations. Non-zero level modules for this Lie algebra were classified in  \cite{[2]}. 
\end{remark}
\vspace{10mm}

\textbf{Acknowledgments}: The authors would like to thank Prof. Eswara Rao for suggesting the problem and some helpful discussion.

\end{document}